# A HYBRID BOUNDARY ELEMENT METHOD FOR ELLIPTIC PROBLEMS WITH SINGULARITIES

George Pashos, Athanasios G. Papathanasiou, Andreas G. Boudouvis

School of Chemical Engineering, National Technical University of Athens, Athens

15780, Greece

**Abstract**

The singularities that arise in elliptic boundary value problems are treated locally by a singular function boundary integral method. This method extracts the leading singular coefficients from a series expansion that describes the local behavior of the singularity. The method is fitted into the framework of the widely used boundary element method (BEM), forming a hybrid technique, with the BEM computing the solution away from the singularity. Results of the hybrid technique are reported for the Motz problem and compared with the results of the standalone BEM and Galerkin/finite element method (GFEM). The comparison is made in terms of the total flux (i.e. the capacitance in the case of electrostatic problems) on the Dirichlet boundary adjacent to the singularity, which is essentially the integral of the normal derivative of the solution. The hybrid method manages to reduce the error in the computed capacitance by a factor of 10, with respect to the BEM and GFEM.

## 1 Introduction

Blamed for serious damages, in many engineering applications, singularities are under extensive computational investigation aiming to explore their origin and predicting their effects. A thorough review of several computational techniques that specialize in the treatment of the singularities of elliptic problems can be found in [1]. The cause of these singularities is found in, either a ‘sudden’ change in the boundary conditions (cf. the Motz problem in section 3), or the existence of a sharp corner in the geometry of the computational domain, both addressed in [1] and similarly treated. As for the abrupt changes in geometry, an example of practical interest would be the electric field singularities of conducting wedges surrounded by dielectrics and vice versa, in electrostatic problems [2]. A case of particular interest is the field singularity that arises in electro-capillary systems at the contact line (abrupt change in geometry), which is investigated in connection with phenomena that limit the electrostatically assisted wetting, such as the contact angle saturation [3].

A large family of techniques treating elliptic problems with singularities accounts for the asymptotic expansion of the singularity. In the present work we restrict ourselves to problems governed by the Laplace equation, which is the most common representative of the family of elliptic boundary value problems. The analysis of the singularity of the Laplace equation posed in a 2-d arbitrary domain – which will be the basis of the proposed method – produces an asymptotic expansion of the solution, $u$, that reads, $u = \sum_{l=1}^{\infty} a_l r^{\mu_l} f_l(\theta)$ where, $\alpha_l$ are the unknown singular coefficients (in fracture mechanics known as generalized stress intensity factors), $r$ is the radial distance from the center of the singularity, $\theta$ is the angle with reference to a boundary of the 2-d domain, $\mu_l$ and $f_l$ are predetermined by the analysis of the singularity (cf. section 2).

The goal of this work is to embed the asymptotic expansion of the solution into a widely used computational method for analyzing physical systems governed by the Laplace equation. The method of choice is the boundary element method (BEM) [4, 5] that is already a very commonly used technique in elasticity and potential problems and it is gaining ground in other types of problems, mainly due to its reduced computational cost, compared with other methods, e.g. the finite element method. Considering that matter, the greatest merit of BEM is that it reduces the dimension of the computational problem by one, i.e. only the boundary of the computational domain is discretized.

The BEM has already been enhanced by singularity techniques; for example in [6-8] the BEM is augmented with singular boundary elements, i.e. elements with special basis functions that account for the singularity instead of the common polynomial basis functions. However, to embed as many leading terms of the asymptotic expansion of the singularity as the basis functions, requires the construction of many-node elements which in turn demands much tedious work.

In this work the BEM is augmented with elements of the technique in [9], referred to as singular function boundary integral method (SFBIM), which employs the asymptotic expansion of $u$ in the solution procedure, with as many singular terms as required. A similar treatment of singularities with a hybrid BEM was presented in [10]. The SFBIM was introduced in [9] and thereafter was efficiently applied to problems governed by the Laplace equation, such as the L-shaped domain problem [11]. Moreover, SFBIM was applied to problems governed by the biharmonic equation, such as the Newtonian stick-slip flow problem [12] and fracture problems with crack singularities [13].

The proposed method is essentially a coupling of BEM and SFBIM that results in a novel hybrid method and will be referred to as hybrid BEM/SFBIM. The effectiveness of the coupling lies in the similarities of the two methods, with both, being boundary integral methods that reduce the dimension of the computational problem. Briefly, the coupling of the two methods is done as follows: The computational domain is decomposed in two subdomains, the first being a small segment of the original domain, surrounding the singular point, where the SFBIM is applied, and the second its complement, where the BEM is applied – the coupling is analyzed in detail in the following sections.

## 2 The hybrid BEM/SFBIM

Consider the Laplace equation posed in a 2-d arbitrary domain, $\Omega$, depicted in Fig. 1. Let $\Gamma_D$ and $\Gamma_N$ be the boundaries on which the Dirichlet and Neumann boundary conditions are imposed, respectively. The boundary of $\Omega$ is smooth except at the points $O_2$, $O_3$, where the boundary forms the angles $\Theta_2$, $\Theta_3$, respectively; $O_1$ lies on a straight boundary segment with $\Theta_1=\pi$, therefore the boundary is smooth at $O_1$. Singularities arise at the points of the boundary at which a sharp corner is formed ($O_2$) or there is a sudden change of the boundary condition ($O_1$) or both ($O_3$). The change of the boundary condition can be between Dirichlet and Neumann or of the same kind, e.g. homogeneous-inhomogeneous Dirichlet. The asymptotic solution near $O_1$, $O_2$ and $O_3$ can be derived through separation of the independent variables, $r$, $\theta$, where $r$ is the radial distance from $O_i$ and $\theta$ is the angle with reference to the boundary. The boundary segments adjacent to $O_1$, $O_2$ and $O_3$ are straight lines in order to provide a simple form of the asymptotic solution. A more general expression for curved boundaries can be found in [1]. In this work, however, we restrict ourselves to the analysis of the simplest cases, that is homogeneous Dirichlet and/or Neumann conditions imposed on straight boundaries adjacent to the singular points. The general expression of the asymptotic solution for Neumann-Dirichlet (with $\theta$ = 0 on the Neumann boundary) conditions is

$$u = u_0(r,\theta) + \sum_{l=1}^{\infty} a_l r^{\mu_l} \cos(\mu_l \theta)$$

and for Dirichlet-Dirichlet conditions

$$u = u_0(r,\theta) + \sum_{l=1}^{\infty} a_l r^{\mu_l} \sin(\mu_l \theta)$$

where $\alpha_l$ are the unknown singular coefficients, $\mu_l$ are the known powers of the singularity that depend on $\Theta_i$ ($\Theta_1$, $\Theta_2$, or $\Theta_3$) and the types of boundary conditions (the details for the extraction of $\mu_l$ are given in [1]), $u_0$ are particular solutions, which vanish for homogeneous boundary conditions, thus giving

$$u = \sum_{l=1}^{\infty} a_l r^{\mu_l} \cos(\mu_l \theta) \tag{1}$$

and

$$u = \sum_{l=1}^{\infty} a_l r^{\mu_l} \sin(\mu_l \theta) \tag{2}$$

with (1) being valid as well for homogeneous Neumann-Neumann conditions.

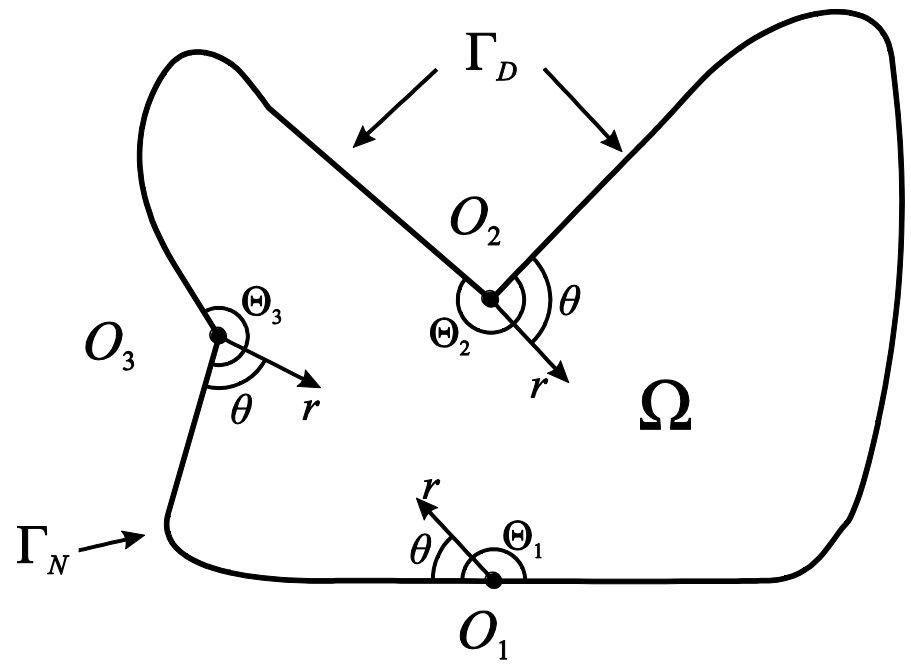

Fig. 1. A 2-d arbitrary domain with Dirichlet ($\Gamma_D$), Neumann ($\Gamma_N$) boundary conditions and singularity source points $O_1$, $O_2$ and $O_3$.

As shown in Fig. 1, multiple singularities can exist in the same domain and we assert ourselves that the local expression of each singularity is valid as it was extracted neglecting the presence of the other singularities. The effect of a pair of singularities on the solution is investigated in section 3.

The original domain of Fig. 1 is decomposed into three non-overlapping subdomains; the small subdomains $\Omega_1$, $\Omega_2$ contain the singularities – for simplicity we neglect the singularity at the point $O_3$ -- and the large subdomain $\Omega_3$ contains the bulk space of $\Omega$ where the effect of the singularity is relatively weak (Fig. 2). The subdomains are separated through a circular segment that surrounds the singular points at a given radius, $R$. The choice of the shape of the artificial inner boundary, that is being circular, is not justified by any other means than simplicity of implementation and the fact that is the most straightforward approach since the strength of the singularity depends on the radial distance from $O_i$. The boundaries of each subdomain are, for $\Omega_1$ $\Gamma = \Gamma_4$, for $\Omega_2$ $\Gamma = \Gamma_5$ and for $\Omega_3$ $\Gamma = \Gamma_1 \cup \Gamma_4 \cup \Gamma_2 \cup \Gamma_5 \cup \Gamma_3$ -- the small straight segments of the boundaries of $\Omega_2$ and $\Omega_3$ are not taken into account for the present analysis due to the properties of the suggested method, as it will be discussed below. The boundaries are grouped with reference to the type of the boundary condition as, $\Gamma_N = \Gamma_1$, $\Gamma_D = \Gamma_2 \cup \Gamma_3$ and the internal boundaries, $\Gamma_4$, $\Gamma_5$, don't fall into either one of the above types of boundaries as they do not have a specified boundary condition. The solution is approximated via (1) and (2) only in $\Omega_1$ and $\Omega_2$, respectively, while in $\Omega_3$ the solution is approximated via standard polynomial interpolation functions that are typically constant or linear.

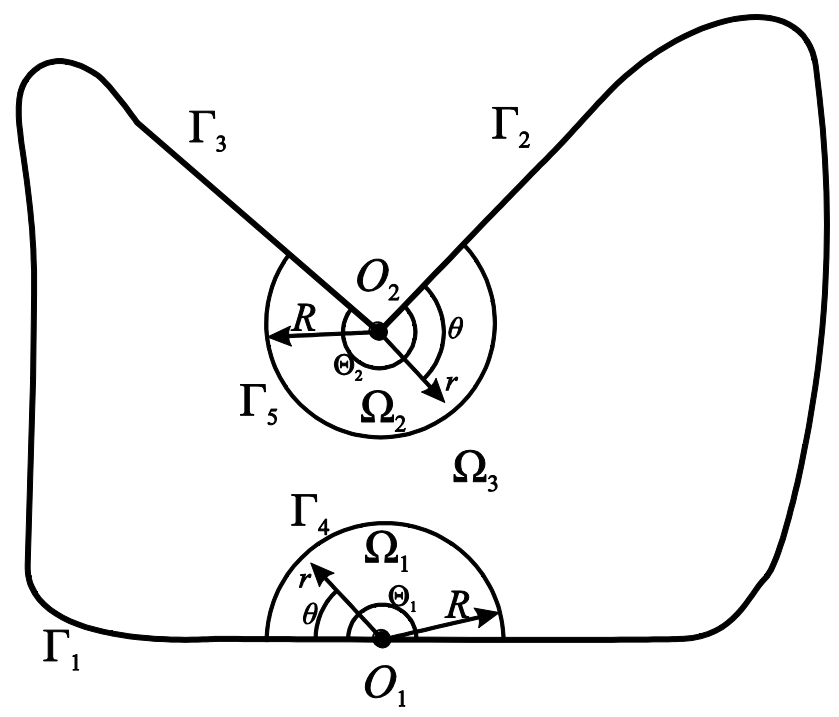

Fig. 2. Domain decomposition for the hybrid boundary integral method with two singularities; $\Gamma_1 = \Gamma_N$, $\Gamma_2 = \Gamma_D$, $\Gamma_3 = \Gamma_D$.

In $\Omega_3$ we apply the standard BEM. In more detail, starting from the fundamental solution of the Laplace equation, the BEM extracts the following boundary integral equation

$$\lambda(\xi,n)u(\xi,n) + \int_\Gamma u(x,y)\frac{\partial G(x,y;\xi,n)}{\partial n}d\Gamma = \int_\Gamma \frac{\partial u(x,y)}{\partial n}G(x,y;\xi,n)d\Gamma \tag{3}$$

where $G(x,y;\xi,n)$ is the free-space Green's function, defined as $G(x,y;\xi,n) = \frac{1}{2\pi}\ln(1\big/\sqrt{(x-\xi)^2+(y-n)^2})$ and $n$ is in the direction of the unit vector normal to the boundary. The parameter $\lambda$ depends on ($\xi$, $n$); if $(\xi, n) \in \Omega_3$ then $\lambda$=1, if $(\xi, n) \notin \Omega_3$ or $\partial\Omega_3$ then $\lambda$=0; if $(\xi, n) \in \partial\Omega_3$ cf. below. The boundaries of $\Omega_3$ are tessellated into a finite number of constant or linear elements, with one node placed in the center of each element or two nodes placed at the endpoints of each element, respectively. The solution, $u(x, y)$, and its derivative, $\partial u(x, y)/\partial n$, on the boundary, $\Gamma$, are approximated in terms of the basis functions, $\Phi^j(x, y)$, as

$$u(x,y) = \sum_{j=1}^{N} u_j \, \Phi^j(x,y), \quad \frac{\partial u(x,y)}{\partial n} = \sum_{j=1}^{N} q_j \Phi^j(x,y) \tag{4}$$

where $u_j$ and $q_j$ are the nodal values of $u$ and its derivative, respectively and $N$ is the total number of nodes. The basis functions, $\Phi^j(x, y)$, are piecewise polynomial functions that can either be constant or linear. In this work we use linear basis functions, unless indicated otherwise.

By collocating the points ($\xi$, $n$) with the nodal positions ($x_i$, $y_i$) and inserting (4) in (3), a discretized version of (4) is derived that reads

$$\begin{aligned} &\lambda(x_i, y_i)u_i + \sum_{j=1}^{N} u_j \int_{\Gamma} \Phi^j(x,y) \frac{\partial G(x,y;x_i,y_i)}{\partial n} d\Gamma = \sum_{j=1}^{N} q_j \int_{\Gamma} \Phi^j(x,y) G(x,y;x_i,y_i) d\Gamma \\ &i = 1, \cdots, N \end{aligned} \tag{5}$$

The set of equations (5) can be written in matrix form

$$\mathbf{Hu} = \mathbf{Gq} \tag{6}$$

where

$$\begin{aligned} H_{ij} &= \int_{\Gamma} \Phi^j(x,y) \frac{\partial G(x,y;x_i,y_i)}{\partial n} d\Gamma, \quad for \quad i \neq j \\ H_{ij} &= \lambda(x_i, y_i) + \int_{\Gamma} \Phi^j(x,y) \frac{\partial G(x,y;x_i,y_i)}{\partial n} d\Gamma, for \quad i = j \\ G_{ij} &= \int_{\Gamma} \Phi^j(x,y) G(x,y;x_i,y_i) d\Gamma \end{aligned} \tag{7}$$

The parameter $\lambda$ in (7) is equal to ½ if the $i$-th node lies in a smooth segment of the boundary, e.g. the middle of an element (that is the case for constant basis functions). When the $i$-th node lies in a corner of the boundary, e.g. the endpoints of the element (that is the case for linear basis functions), $\lambda$ depends on the angle formed by the elements that are adjacent to the $i$-th node. It is convenient to circumvent the explicit computation of $\lambda$ with a simple technique without any loss of generality of the method [4]. This is done by applying a uniform $u$ along the boundary, which in turn bounds the normal derivative of the solution, $q$, to be zero. Thus, (6) becomes

$$\mathbf{Hu} = \mathbf{0}$$

which in turn provides the diagonal elements of **H** using the computed off-diagonal elements,

$$H_{ii} = -\sum_{\substack{j=1 \\ j \neq i}}^{N} H_{ij}$$

instead of using the second equation of (7) that requires the computation of $\lambda$.

On each boundary of the domain there is a specified boundary condition, either Dirichlet where the $u_j$ are defined or Neumann where the $q_j$ are defined, excluding the internal boundaries where no boundary condition is applied. Thus, we can reorder (6) so that the unknowns lie in the left hand side of the equation

$$\mathbf{Ax} = \mathbf{b} \tag{8}$$

Based on the setup of the problem of Fig. 2, the matrix and the vectors of the unknowns are, respectively

$$\begin{gathered} \mathbf{A} = [\ [\mathbf{H}]_{j\in\Gamma_1} \quad -[\mathbf{G}]_{j\in\Gamma_2} \quad -[\mathbf{G}]_{j\in\Gamma_3} \quad [\mathbf{H}]_{j\in\Gamma_4} \quad [\mathbf{H}]_{j\in\Gamma_5} \quad -[\mathbf{G}]_{j\in\Gamma_4} \quad -[\mathbf{G}]_{j\in\Gamma_5}\ ] \\ \mathbf{x} = \left[\mathbf{u}_{\Gamma_1}, \mathbf{q}_{\Gamma_2}, \mathbf{q}_{\Gamma_3}, \mathbf{u}_{\Gamma_4}, \mathbf{u}_{\Gamma_5}, \mathbf{q}_{\Gamma_4}, \mathbf{q}_{\Gamma_5}\right]^{\mathrm{T}}. \end{gathered} \tag{9}$$

**A** is a $N \times (N + M_{\Gamma_4} + M_{\Gamma_5})$ matrix, where $M_{\Gamma_4}$ and $M_{\Gamma_5}$ are the number of nodes on $\Gamma_4$ and $\Gamma_5$, respectively. The extra $M_{\Gamma_4} + M_{\Gamma_5}$ columns correspond to the two rightmost set of columns of **A**, $-[\mathbf{G}]_{j\in\Gamma_4}$, $-[\mathbf{G}]_{j\in\Gamma_5}$, which are gathered in the LHS of (8) since $\Gamma_4$ and $\Gamma_5$ are internal boundaries and thus, both **u** and **q** are unknown. The extra $M_{\Gamma_4} + M_{\Gamma_5}$ equations needed, will be provided by the coupling with the SFBIM that is applied on the subdomains, $\Omega_1$ and $\Omega_2$.

In $\Omega_1$ and $\Omega_2$ the solution, $u$, is approximated by (1) or (2) that can be rewritten

$$u = \sum_{l=1}^{\infty} a_l W^l(r,\theta) \tag{10}$$

where $W^l$ are harmonic functions of $r$ and $\theta$. The SFBIM incorporates the $N_\alpha$ leading terms of (10), however, the proposed method performs sufficiently well with just a few leading terms as it will be seen in section 3. Thus, (10) is rewritten

$$u = \sum_{l=1}^{N_a} a_l W^l(r,\theta) \tag{11}$$

The derivative of $u$ normal to the boundaries of $\Omega_1$ and $\Omega_2$ is approximated by

$$\frac{\partial u}{\partial n} = \sum_{l=1}^{N_a} a_l \frac{\partial W^l(r,\theta)}{\partial n} \tag{12}$$

The Laplace equation is weighted with $W^k$ and Green's theorem is applied twice. The double integral that contains the term, $\nabla^2 W^k$, vanishes since $W^k$ are harmonic and thus, the following boundary integral equation is obtained

$$\int_\Gamma \frac{\partial u}{\partial n} W^k d\Gamma - \int_\Gamma u \frac{\partial W^k}{\partial n} d\Gamma = 0 \qquad k = 1,\ldots,N_a \tag{13}$$

Next, we apply (13) on the subdomains $\Omega_1$ and $\Omega_2$. The solution, $u$, is approximated on boundaries of the subdomains via (11). Its derivative, however, is approximated on the straight segments of $\Omega_1$ and $\Omega_2$ via (12) and on the inner artificial boundaries, $\Gamma_4$ and $\Gamma_5$, via polynomial basis functions (cf. Eq. (4)). This approach provides the equality constraints for the derivative of the solution at the artificial inner boundaries ($C^1$ continuity constraints). The resulting set of integral equations from (13) reads

$$\begin{gathered}\sum_{j=1}^{N} q_j \int_{\Gamma\subset\partial\Omega_3} \Phi^j W^k d\Gamma + \sum_{l=1}^{N_a} a_l \int_{\Gamma\not\subset\partial\Omega_3} \frac{\partial W^l}{\partial n} W^k d\Gamma - \sum_{l=1}^{N_a} a_l \int_\Gamma W^l \frac{\partial W^k}{\partial n} d\Gamma = 0 \\ k = 1,\ldots,N_a\end{gathered} \tag{14}$$

The second boundary integral of (14) vanishes because either $W$ or $\partial W/\partial n$ is zero – depending on the boundary condition (homogeneous Neumann, $\partial W/\partial n|_{\Gamma_N} = 0$ and homogeneous Dirichlet, $W|_{\Gamma_D} = 0$) -- on the straight boundary segments of $\Omega_1$ and $\Omega_2$; the same applies for the third boundary integral of (14) for $\Gamma \not\subset \partial\Omega_3$. Thus, (14) becomes

$$\sum_{j=1}^{N} q_j \int_{\Gamma\subset\partial\Omega_3} \Phi^j W^k d\Gamma - \sum_{l=1}^{N_a} a_l \int_{\Gamma\subset\partial\Omega_3} W^l \frac{\partial W^k}{\partial n} d\Gamma = 0 \qquad k = 1,\ldots,N_a \tag{15}$$

Here are introduced along with the $N_\alpha$ equations of (15), also $N_\alpha$ unknowns (the leading singular coefficients) for each subdomain that contains a singularity; for the subdomains $\Omega_1$ and $\Omega_2$ are introduced $N_{a,\Omega_1}$ and $N_{a,\Omega_2}$ unknowns, respectively. The rest of the constraints are provided by the matching requirement, that is to equalize the approximations of $u$, weighted by the basis functions $\Phi^j$, of the BEM and SFBIM along the boundaries, $\Gamma_4$ and $\Gamma_5$ ($C^0$ continuity constraints)

$$\sum_{j=1}^{N} u_j \int_{\Gamma \subset \partial\Omega_3} \Phi^j \Phi^i d\Gamma - \sum_{l=1}^{N_a} a_l \int_{\Gamma \subset \partial\Omega_3} W^l \Phi^i d\Gamma = 0 \qquad i = 1,\ldots,M \tag{16}$$

where $M = M_{\Gamma_4}, M_{\Gamma_5}$ and $M_{\Gamma_4}, M_{\Gamma_5}$ are the number of elements on $\Gamma_4$ and $\Gamma_5$, respectively. From (16), $M_{\Gamma_4} + M_{\Gamma_5}$ constraints are introduced to the problem, and overall, from (5), (15) and (16) we gather $N + M_{\Gamma_4} + M_{\Gamma_5} + N_{a,\Omega_1} + N_{a,\Omega_2}$ equations with the same number of unknowns. The system (8) is completed with (15) and (16) and **A**, **x** are in expanded form (cf. (9) for the corresponding incomplete system)

$$\mathbf{A} = \begin{bmatrix} [\mathbf{H}]_{j\in\Gamma_1} & -[\mathbf{G}]_{j\in\Gamma_2} & -[\mathbf{G}]_{j\in\Gamma_3} & [\mathbf{H}]_{j\in\Gamma_4} & [\mathbf{H}]_{j\in\Gamma_5} & -[\mathbf{G}]_{j\in\Gamma_4} & -[\mathbf{G}]_{j\in\Gamma_5} & \mathbf{0} & \mathbf{0} \\ \mathbf{0} & \mathbf{0} & \mathbf{0} & \mathbf{0} & \mathbf{0} & \int_{\Gamma_4} \mathbf{\Phi W} d\Gamma & \mathbf{0} & -\int_{\Gamma_4} \mathbf{W} \frac{\partial \mathbf{W}}{\partial n} d\Gamma & \mathbf{0} \\ \mathbf{0} & \mathbf{0} & \mathbf{0} & \mathbf{0} & \mathbf{0} & \mathbf{0} & \int_{\Gamma_5} \mathbf{\Phi W} d\Gamma & \mathbf{0} & -\int_{\Gamma_5} \mathbf{W} \frac{\partial \mathbf{W}}{\partial n} d\Gamma \\ \mathbf{0} & \mathbf{0} & \mathbf{0} & \int_{\Gamma_4} \mathbf{\Phi\Phi} d\Gamma & \mathbf{0} & \mathbf{0} & \mathbf{0} & -\int_{\Gamma_4} \mathbf{W\Phi} d\Gamma & \mathbf{0} \\ \mathbf{0} & \mathbf{0} & \mathbf{0} & \mathbf{0} & \int_{\Gamma_5} \mathbf{\Phi\Phi} d\Gamma & \mathbf{0} & \mathbf{0} & \mathbf{0} & -\int_{\Gamma_5} \mathbf{W\Phi} d\Gamma \end{bmatrix}$$

$$\mathbf{x} = \left[\mathbf{u}_{\Gamma_1}, \mathbf{q}_{\Gamma_2}, \mathbf{q}_{\Gamma_3}, \mathbf{u}_{\Gamma_4}, \mathbf{u}_{\Gamma_5}, \mathbf{q}_{\Gamma_4}, \mathbf{q}_{\Gamma_5}, \boldsymbol{\alpha}_{\Omega_1}, \boldsymbol{\alpha}_{\Omega_2}\right]^{\mathrm{T}}$$

where $\boldsymbol{\alpha}_{\Omega_1}$ and $\boldsymbol{\alpha}_{\Omega_2}$ are the vectors of the leading singular coefficients of the singularities that are contained in $\Omega_1$ and $\Omega_2$, respectively, with $\boldsymbol{\alpha}_{\Omega_1} = [a_1, \ldots, a_{N_{a,\Omega_1}}]$ and $\boldsymbol{\alpha}_{\Omega_2} = [a_1, \ldots, a_{N_{a,\Omega_2}}]$.

## 3 Numerical experiments

The proposed technique is applied on the standard benchmark Motz problem [14, 15] (Fig. 3). The problem is governed by the Laplace equation posed in a rectangular domain, $\Omega = [-1,1] \times [0,1]$ that is divided into five boundaries, denoted $\Gamma_1,\ldots, \Gamma_5$. The singularity is centered at the origin of the axes at the intersection of the boundaries $\Gamma_1$, $\Gamma_5$, where there is a sudden change in the boundary condition from $u|_{\Gamma_5} = 0$ to $\partial u / \partial n|_{\Gamma_1} = 0$. On $\Gamma_3$, $\Gamma_4$, is applied the homogeneous Neumann boundary condition, $\partial u / \partial n|_{\Gamma_3 \cup \Gamma_4} = 0$ and on $\Gamma_2$ a Dirichlet boundary condition, $u|_{\Gamma_2} = 500$.

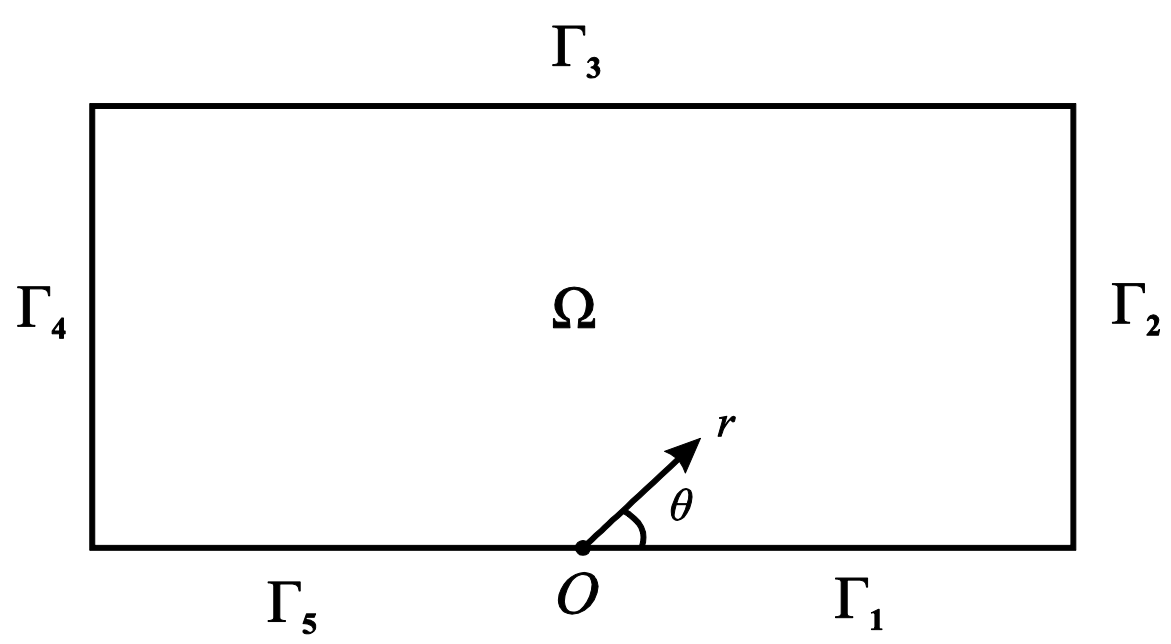

Fig. 3. The Motz problem domain; $\Omega = [-1,1] \times [0,1]$.

The asymptotic solution of the singularity for the Motz problem is given by the infinite series

$$u = \sum_{l=1}^{\infty} a_l r^{\frac{2l-1}{2}} \cos\left[\left(\frac{2l-1}{2}\right)\theta\right] \tag{17}$$

which is exact for the entire domain and thus, true for any subdomain that includes $O$. Therefore, $\Omega$ is decomposed into two non-overlapping subdomains as shown in Fig. 4 where $\Omega_1$ contains the singularity and $\Omega_2$ contains the bulk space of the original domain. The harmonic functions for the Motz problem for the subdomain $\Omega_2$ are

$$W^k = r^{\frac{2k-1}{2}} \cos\left[\left(\frac{2k-1}{2}\right)\theta\right] \tag{18}$$

which are valid for the entire domain; however, that is not a requirement for the proposed method. Only the leading $N_\alpha$ harmonic functions are incorporated in the solution procedure, i.e. $k = 1,\ldots,N_a$, where $N_\alpha$ is the number of the leading singular coefficients.

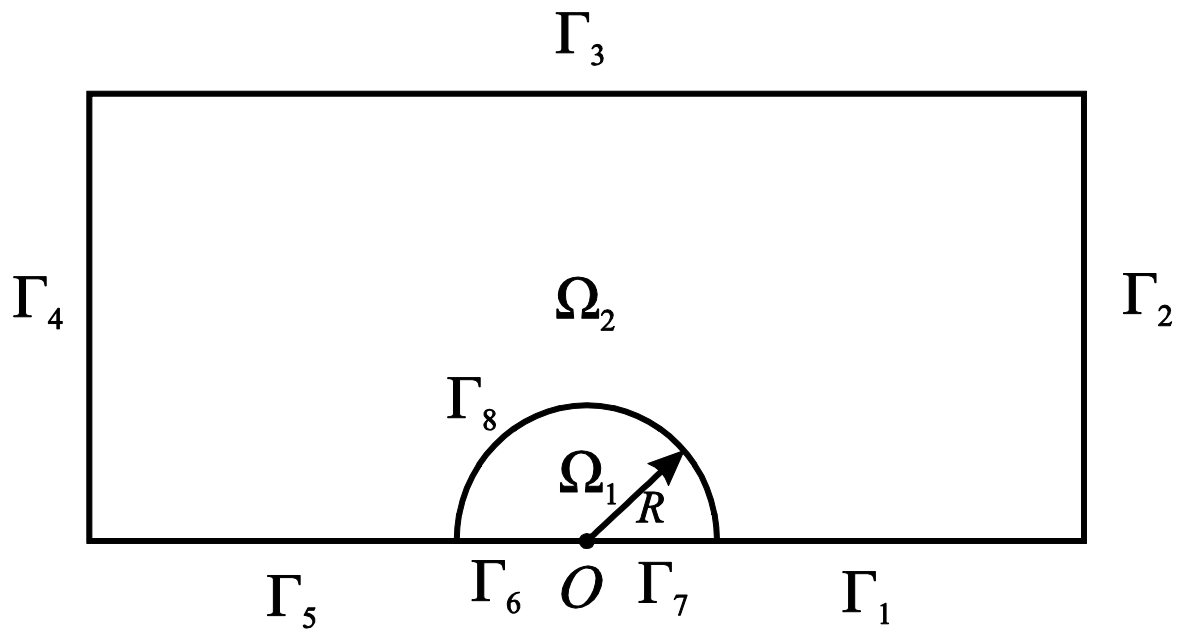

Fig. 4. Domain decomposition of the Motz problem; $\Omega \equiv \Omega_1 \cup \Omega_2 = [-1,1] \times [0,1]$.

The application of the BEM on the domain $\Omega_2$ depicted in Fig. 4 produces a system of equations as in (8) for a problem with arbitrary domain. The LHS reads

$$\mathbf{A} = [\ [\mathbf{H}]_{j\in\Gamma_1} \quad -[\mathbf{G}]_{j\in\Gamma_2} \quad [\mathbf{H}]_{j\in\Gamma_3} \quad [\mathbf{H}]_{j\in\Gamma_4} \quad -[\mathbf{G}]_{j\in\Gamma_5} \quad [\mathbf{H}]_{j\in\Gamma_8} \quad -[\mathbf{G}]_{j\in\Gamma_8}\ ] \tag{19}$$

The system is augmented with $N_a + M$ equations resulting from the application of the continuity constraints on $\Omega_1$ of Fig. 4

$$\sum_{j=1}^{N} q_j \int_{\Gamma_8} \Phi^j W^k d\Gamma - \sum_{l=1}^{N_a} a_l \int_{\Gamma_8} W^l \frac{\partial W^k}{\partial n} d\Gamma = 0 \qquad k = 1,\ldots,N_a \tag{20}$$

$$\sum_{j=1}^{N} u_j \int_{\Gamma_8} \Phi^j \Phi^i d\Gamma - \sum_{l=1}^{N_a} a_l \int_{\Gamma_8} W^l \Phi^i d\Gamma = 0 \qquad i = 1,\ldots,M \tag{21}$$

The resulting system from (19), (20) and (21) has a size $N + M + N_a$ vector of unknowns that reads $\mathbf{x} = [\mathbf{u}_{\Gamma_1}, \mathbf{q}_{\Gamma_2}, \mathbf{u}_{\Gamma_3}, \mathbf{u}_{\Gamma_4}, \mathbf{q}_{\Gamma_5}, \mathbf{u}_{\Gamma_8}, \mathbf{q}_{\Gamma_8}, \boldsymbol{\alpha}]^{\mathrm{T}}$.

The efficiency of the proposed method is evaluated in terms of the computed singular coefficients, $\alpha_l$, compared with their exact values found in [9]. Table 1 shows a set of results for the leading singular coefficients with typical solution parameters such as the discretization of the boundaries and the radius of $\Omega_1$. The internal artificial boundary and the external boundaries of the Motz problem are discretized uniformly, e.g. if $N = 100$ then each side of the rectangular domain has 50 equally sized elements and if $M = 10$ then the boundary $\Gamma_8$ is represented by an even polygon with 10 sides; if linear or constant elements are incorporated the boundaries are represented as straight segments even though they may be curved. The results of the proposed method are in good agreement with the exact values of the two leading singular coefficients beyond which the discrepancies become so large that can possibly exceed several orders of magnitude – it should be noted that large deviations of the computed singular coefficients, from the third and above, do not affect the computed $u_j$ or $q_j$. In addition, regarding the solution in the subdomain where the SFBIM is applied, there is no 'noticeable' difference of the solution (that has the form of (11)) when the third leading singular coefficient and above is largely miscalculated – this is restricted to relatively small values of $R$ with respect to the size of the computational domain, for this case ~0.1. To further illustrate this point, we can apply an arbitrary value of $\alpha_3$, $\alpha_4$, etc. and the deviation of computed solution (from (11)) from the exact (again from (11) using the exact $\alpha_l$) will not exceed margins of practical interest; this is the case when $R$ is relatively small. When $R$ is further increased, $N_\alpha$ has to be increased as well, in order to achieve the same levels of accuracy (cf. below for a quantitative investigation (Fig. 5)). In problems with practical interest, a sensible choice of $R$ requires only $N_\alpha$=1 or 2, since its value would be relatively small. However, to evaluate the performance of the proposed method $N_\alpha$ is increased up to 10, while maintaining a small value of $R$. Further increase of $N_\alpha$ would be unnecessary since it has a little effect on the computed $\alpha_l$ as seen in Table 2. This increase is limited however to a relatively small range of values of $N_\alpha$ because otherwise the matrix **A** becomes ill-conditioned and the accuracy of the solution is compromised.

Table 1

Five leading singular coefficients of the Motz problem derived with the proposed method and their exact values; $N$ = 500 (total number of elements), $M$ = 100 (number of elements on inner artificial boundary), $R$ =0.1, $N_\alpha$ = 5.

| | $\alpha_1$ | $\alpha_2$ | $\alpha_3$ | $\alpha_4$ | $\alpha_5$ |
|---|---|---|---|---|---|
| Values of the proposed method | 401.067 | 84.7409 | 32.5103 | -153.517 | 994.894 |
| Exact value | 401.162 | 87.6559 | 17.2379 | -8.07121 | 1.44027 |

Table 2

Dependence of the $\alpha_l$ on $N_\alpha$; $N$ = 100, $M$ = 20, $R$ = 0.1.

| | $N_\alpha$=1 | $N_\alpha$=2 | $N_\alpha$=5 | $N_\alpha$=10 | Exact value |
|---|---|---|---|---|---|
| $\alpha_1$ | 401.240 | 401.223 | 401.231 | 401.221 | 401.162 |
| $\alpha_2$ | - | 82.9425 | 82.9647 | 83.0262 | 87.6559 |
| $\alpha_3$ | - | - | 39.9397 | 39.4761 | 17.2379 |
| $\alpha_4$ | - | - | -213.750 | -210.720 | -8.07121 |
| $\alpha_5$ | - | - | 1940.94 | 1915.81 | 1.44027 |
| $\alpha_6$ | - | - | - | -11668.5 | 0.331054 |
| $\alpha_7$ | - | - | - | 158037 | 0.275437 |
| $\alpha_8$ | - | - | - | -806343 | -0.0869329 |
| $\alpha_9$ | - | - | - | 13477300 | 0.0336048 |
| $\alpha_{10}$ | - | - | - | -65620700 | 0.0153843 |

The findings above suggest that a proper computational practice regarding $N_\alpha$ is to incorporate only the leading singular coefficients that are needed (see below). Exceeding a reasonably small value of $N_\alpha$ does not improve the solution, although it causes larger computational cost by further ill-conditioning the matrix **A**.

The required value of $N_\alpha$ depends on the quantity we wish to compute with high accuracy and also the radius, $R$. A very useful quantity for many engineering applications is the total flux on Dirichlet boundaries. For example, in electrostatics, the capacitance, which is traditionally used instead of the term total flux, is a quantity of primary interest [16, 17]. In this work and throughout the text we adopt the term capacitance (instead of total flux), denoted $C$, computed on boundaries with Dirichlet boundary conditions, e.g. $\Gamma_6$ of Fig. 4. It is defined as

$$C = \int_{\Gamma_D} \frac{\partial u}{\partial n} d\Gamma \tag{22}$$

If the boundary belongs to a subdomain that contains a singularity and specifically for the subdomain $\Omega_1$ of the Motz problem of Fig. 4, then $\partial u / \partial n$ is approximated by (12) and (18). Thus, (22) on $\Gamma_6$ gives

$$C = \sum_{l=1}^{N_\alpha} a_l K_l \tag{23}$$

With

$$K_l = \sin\left(\frac{2l-1}{2}\pi\right) R^{\frac{2l-1}{2}} \tag{24}$$

where $K_l$ can be viewed as a weighting term for the contribution of $\alpha_l$ in the computation of $C$. In the case that $R$=1, then $|K_l| = 1$ and the contribution of each $\alpha_l$ is equal and therefore, to compute $C$ with accuracy of three significant digits requires at least the five leading $\alpha_l$ as seen in Table 1 (the exact values); this can be seen by summing the exact values of Table 1, multiplied with the appropriate $K_l$ ($K_l$ =1 or -1). However, this is the extreme case where $R$ is equal to the size of the original computational domain. In the usual range of values of $R$ the required $N_\alpha$ for the accurate computation of $C$ is restricted to less than five $\alpha_l$ as it will be seen below.

Table 3
Dependence of $K_l$ on $R$.

| | $K_1$ | $K_2$ | $K_3$ | $K_4$ | $K_5$ |
|---|---|---|---|---|---|
| $R$=0.9 | 0.9486 | -0.8538 | 0.7684 | -0.6915 | 0.6224 |
| $R$=0.8 | 0.8944 | -0.7155 | 0.5724 | -0.4579 | 0.3663 |
| $R$=0.7 | 0.8366 | -0.5856 | 0.4099 | -0.2869 | 0.2008 |
| $R$=0.6 | 0.7745 | -0.4647 | 0.2788 | -0.1673 | 0.1003 |
| $R$=0.5 | 0.7071 | -0.3535 | 0.1767 | -0.08838 | 0.04419 |
| $R$=0.4 | 0.6324 | -0.2529 | 0.1011 | -0.04047 | 0.01619 |
| $R$=0.3 | 0.5477 | -0.1643 | 0.04929 | -0.01478 | 0.004436 |
| $R$=0.2 | 0.4472 | -0.08944 | 0.01788 | -0.003577 | 0.0007155 |
| $R$=0.1 | 0.3162 | -0.03162 | 0.003162 | -0.0003162 | 0.00003162 |

The $K_l$ values for different values of $R$ (Table 3) are derived from the initially known form of the solution and can be used as a 'loose' criterion for the number of $\alpha_l$ needed for the good approximation of $C$ with respect to $R$. For example, for $R = 0.3$ a choice of $N_\alpha$ with relatively good accuracy would be $N_\alpha$=2, given that $K_3 = 0.04929$ compared to $K_1 = 0.5477$. The above criterion is only indicative since the value of $\alpha_l$ that multiplies $K_l$ is unknown. However, it should also be taken into account that the absolute values of $\alpha_l$ decrease when $l$ increases (for the Motz problem). For example, even though $K_1 = 0.94$ and $K_5 = 0.62$ for $R = 0.9$ are very close, the contribution in $C$ of the first term in (23) is ~380 while the contribution of the fifth term is ~0.89. In realistic applications $R$ should only be a small portion of the size of the original computational domain. However, for extensively analyzing the proposed method, $R$ is increased up to 0.9.

The exact dependence of $C$ on $R$ for the Motz problem and with different values of $N_\alpha$ is shown in Fig. 5. The vertical axis corresponds to the relative error of the computed capacitance with respect to the exact for the given $R$, defined as

$$E = \left|C_{ex} - C\right| / C_{ex} \times 100\% \tag{25}$$

where $C_{ex}$ is computed by (23) with the ten exact leading singular coefficients of the Motz problem (cf. Table 2). All of the curves that correspond to different values of $N_\alpha$ follow the same trend, declining at small values of $R$ then reaching a minimum and start climbing again, some of which reaching error values that exceed 10% -- for $N_\alpha$ =1 we don't observe the above behavior, or at least for the range of $R$ that was examined (0.01-0.9). Examining the graph of Fig. 5 we can extract the optimal values of $N_\alpha$ for intervals of $R$ that minimize $E$. For increasing $N_\alpha$ with increments of 1 from $N_\alpha$=1 to $N_\alpha$=5 the optimal intervals of $R$ in the same order are: $R$=0-0.05, $R$=0.05-0.15, $R$=0.15-0.25, $R$=0.25-0.4, $R$=0.4-.

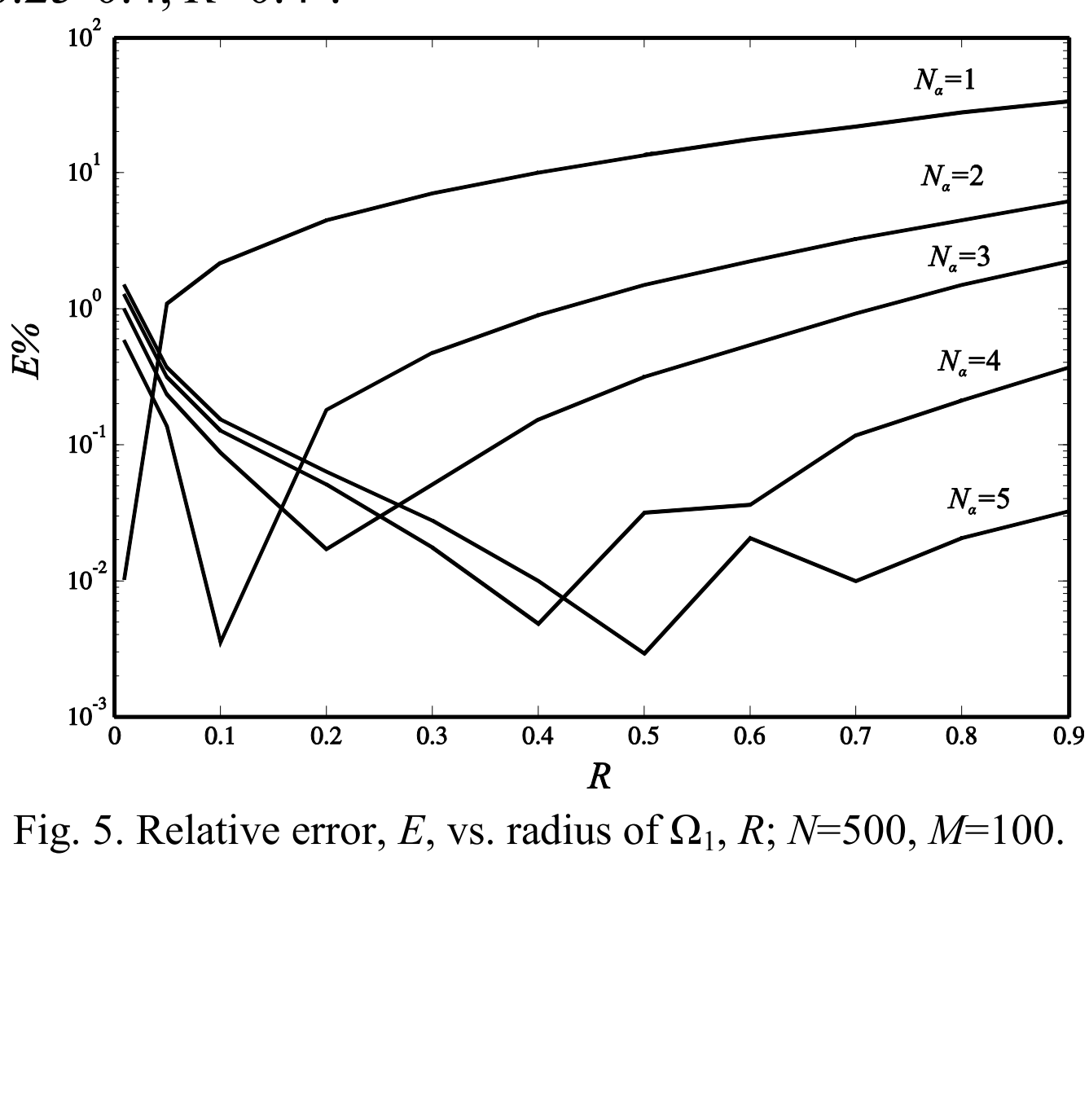

Fig. 5. Relative error, $E$, vs. radius of $\Omega_1$, $R$; $N$=500, $M$=100.

The proposed method is tested against the standalone BEM and the Galerkin/finite element method (GFEM) [18, 19] in terms of the computed capacitance of the Dirichlet boundary, $\Gamma_5 \cup \Gamma_6$ of Fig. 4. On $\Gamma_6$, $C_{\Gamma_6}$ is computed from (23) and on $\Gamma_5$ since the derivative of the solution is approximated via (4) (linear basis functions), a standard trapezoidal integration is adequate for the computation of $C_{\Gamma_5}$. The total capacitance is then $C = C_{\Gamma_5} + C_{\Gamma_6}$. For the standard BEM and GFEM only the trapezoidal rule is applied on the whole lower left boundary. The computed $C$ is compared with the exact $C$ on $\Gamma_5 \cup \Gamma_6$ computed from (23) with $N_\alpha$=10. The relative error is then computed from (25), with $C_{ex}$=340.30. For the various methods to level with each other, we chose a uniform discretization on the boundaries for the BEM and the proposed method, and uniform structured mesh for the GFEM with equally discretized boundaries; e.g. a 40×40 uniform GFEM mesh corresponds to 40 elements per side of Ω that corresponds to $N$ = 160 for the BEM. The current GFEM employs bilinear basis functions and thus, the 40x40 mesh gives 1681 degrees of freedom (DOF), while the BEM gives 160, yet for both of them the number of elements on $\Gamma_5$ of Fig. 4, denoted $N_{\Gamma_5}$ is equal to 20; $\Gamma_6$ does not exist when the BEM or GFEM is applied. The difference between GFEM and BEM in terms of $C$ on $\Gamma_5$ is relatively small, with $E$ following the same pattern with respect to $N_{\Gamma_5}$ (Fig. 6). The GFEM manages to attain a value of $E$ ~3.2% for $N_{\Gamma_5} = 300$ that corresponds to DOF=361201 while the BEM attains the same value for $N_{\Gamma_5} = 400$ and DOF=3200 that also corresponds to $N$=3200; however, the increased DOF of the GFEM is compensated by the sparsity of its matrices, while their counterparts of the BEM are dense. For the comparison of BEM and GFEM with the proposed hybrid method, three different values of $R$ are used that are accompanied by the optimal $N_\alpha$ that minimize $E$ (cf. Fig. 5). Moreover, as in the BEM method the hybrid method uses, $N = 4N_{\Gamma_5 \cup \Gamma_1} + M = 8N_{\Gamma_5} + M$ and for further simplicity, $M$ is constrained as, $M = 2N_{\Gamma_5}$ and thus, $\text{DOF} = N + M + N_a = 12N_{\Gamma_5} + N_a$. For all values of $R$ the proposed method even for small $N_{\Gamma_5}$ achieves small values of $E$, ~1-2%. This is in contrast with the results from Fig. 5 (which are even smaller), due to $C$ in Fig. 5 being computed only on $\Gamma_6$ while in Fig. 6, $C$ is computed on $\Gamma_5 \cup \Gamma_6$.

We can safely assume that the discrepancies between $C_{\Gamma_6}$ and $C_{\Gamma_5 \cup \Gamma_6}$ are related to the integration on $\Gamma_6$, i.e. the boundary treated by the BEM. Moreover, the decrease of $E$ of the proposed method with respect to the increasing $R$ is due to the decrease of $\Gamma_5$ and thus, the contribution of $\Gamma_5$ on the total capacitance. The miscalculation of the integral on $\Gamma_5$ is largely affected by the oscillations of the solution $q_j$ close to the ends of $\Gamma_5$ (Fig. 7), which is a common behavior of the BEM close to corners of domains or at the points where the boundary conditions change. To further illustrate the effect of the oscillations of $q$ on $C_{\Gamma_5}$ and thus, on the total $C$, we seek to eliminate them by simply excluding the solution near the ends of $\Gamma_5$ and fitting a curve based on a nonlinear regression model applied on the remaining values of $q_j$ (Fig. 7). The integral of the fitted curve is then the corrected $C_{\Gamma_5}$. This is done at the post-processing stage and does not belong to the core of the proposed method, but is done only to justify the above remark. However, there are techniques for the BEM, e.g. discontinuous elements at the corners, which provide more accurate solution near the corners, but are not incorporated in the present work.

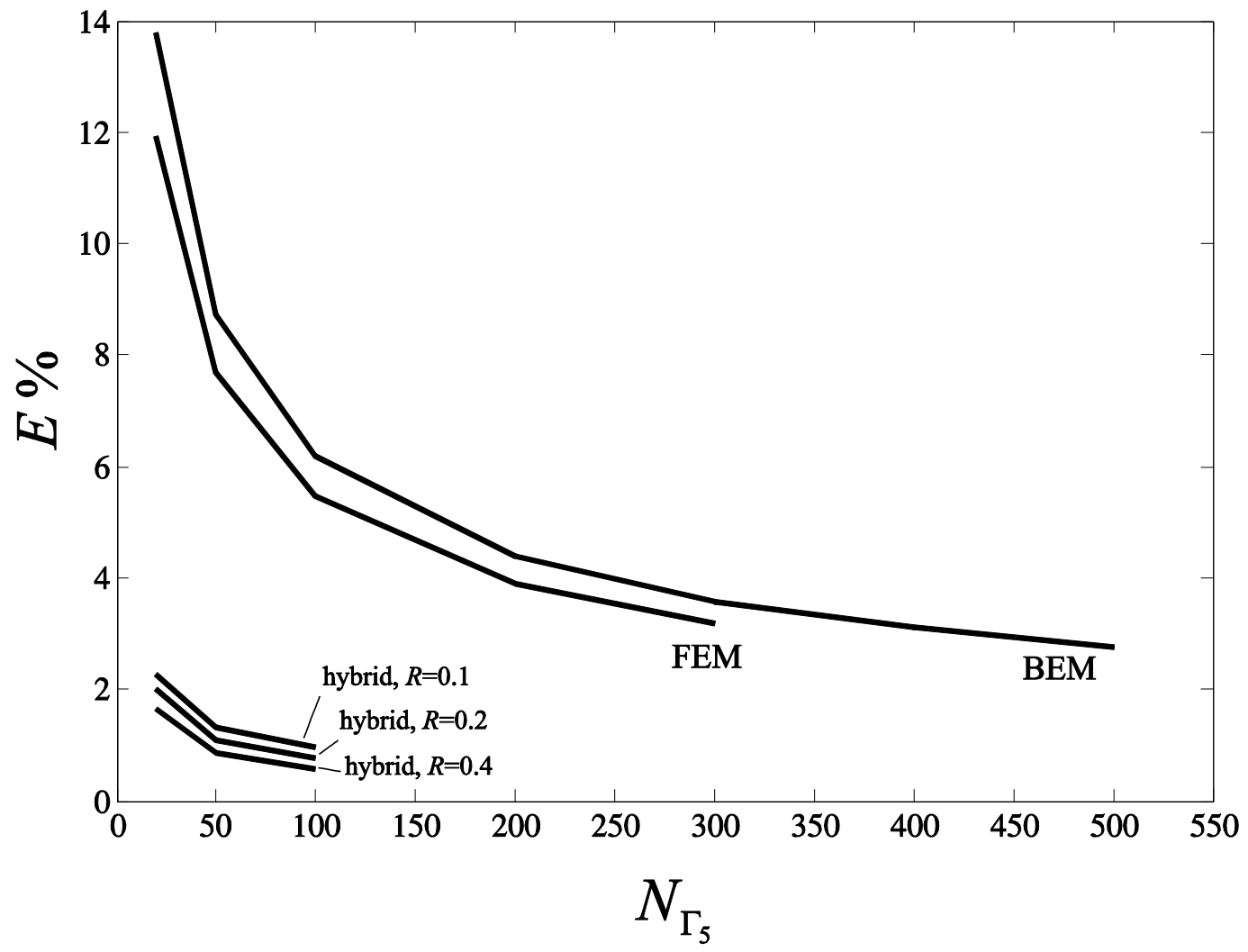

Fig. 6. Relative error, $E$, vs. $N_{\Gamma_5}$ ; comparison between the GFEM, BEM and hybrid integral method.

The relative error, $E$, based on corrected capacitance, $C_{\Gamma_5}$ , which in turn is computed from the integral of the fitted curve of Fig. 7, is presented in Fig. 8. There is a clear reduction of $E$ to levels that are seen in Fig. 5. This should provide sufficient proof that $E$ computed on $\Gamma_5 \cup \Gamma_6$ is relatively large due to the oscillations of the solution, $q_j$, on $\Gamma_5$ near its corners.

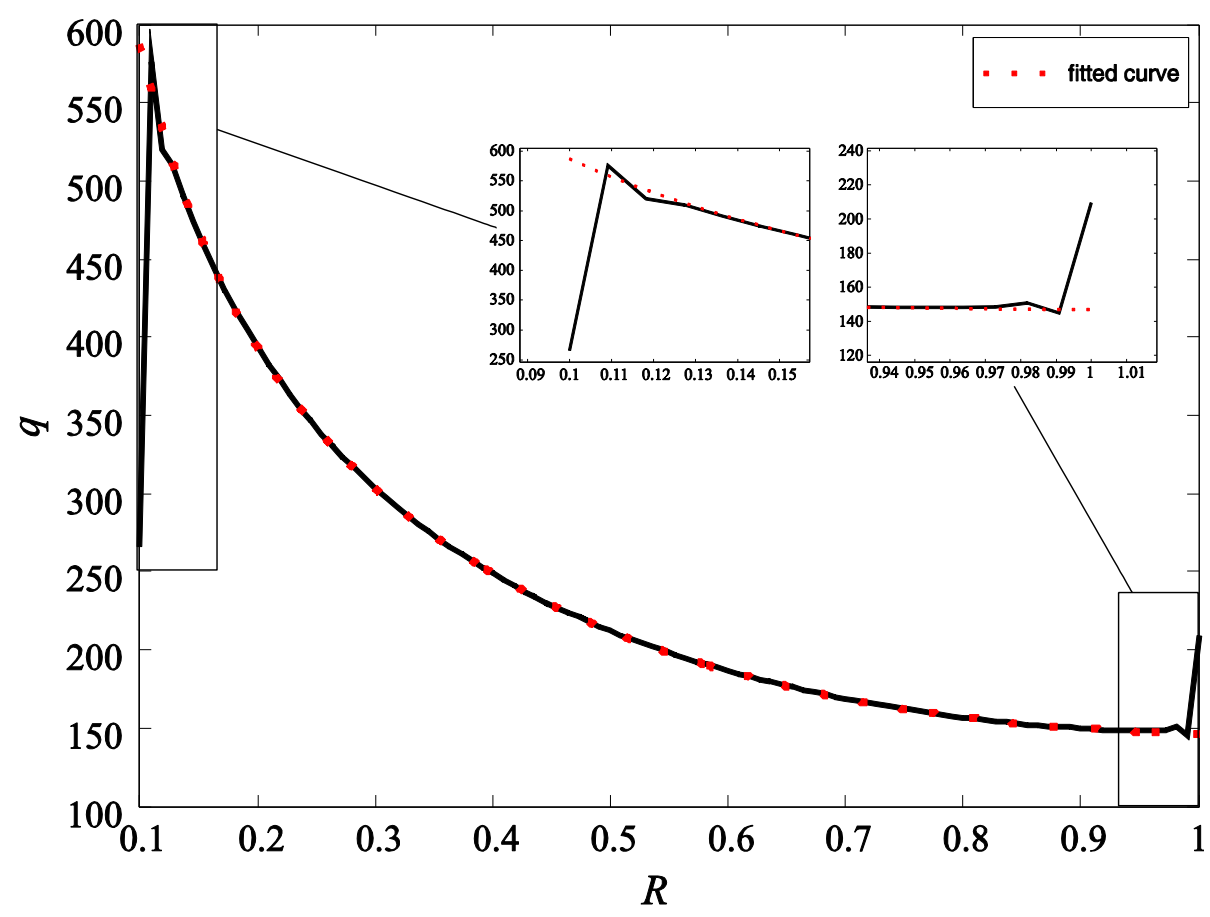

Fig. 7. Solution $q_j$ of the proposed method on $\Gamma_5$ (solid line) and fitted curve (dots); $N_{\Gamma_5}$ =100.

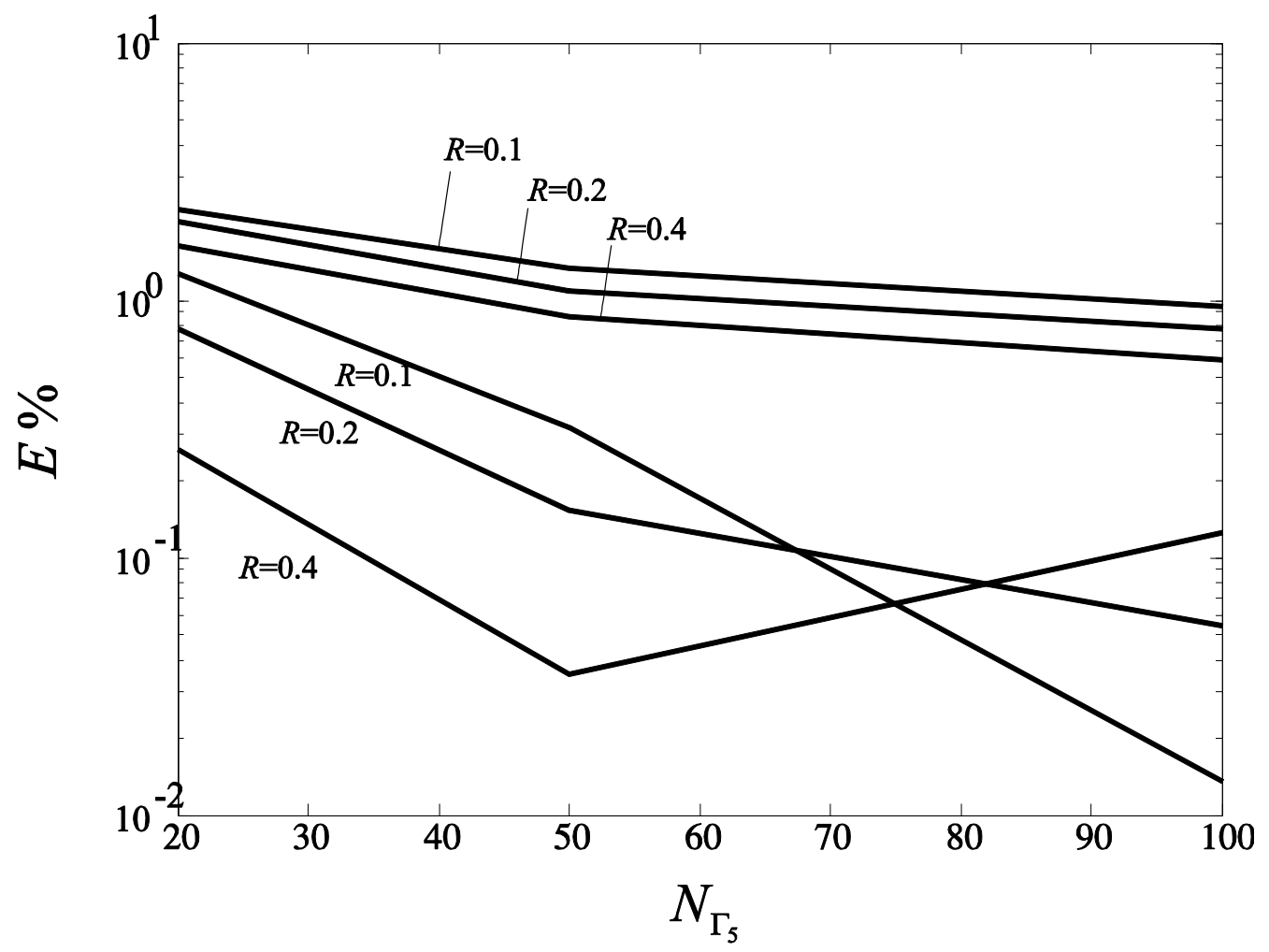

Fig. 8. Relative error, $E$, vs. $N_{\Gamma_5}$ ; solid lines are identical to the ones in Fig. 6 (hybrid method) and dotted lines are derived from (25), implementing the corrected $C$ that is computed by the fitting of $q_j$ on $\Gamma_5$.

The proposed method is tested on a computational domain with two singularities. The problem is again the Motz problem even though it contains only one prominent singularity. However, every corner of the computational domain of the Motz problem can potentially impose a singularity that varies in strength, but comparatively to the singularity at (0, 0) is negligible. Despite this, we treat the area at the upper left corner as we treated the area at (0, 0). The computational domain is decomposed as in Fig. 9. Homogeneous Neumann-Neumann conditions are imposed on the straight segments of the boundary of $\Omega_2$ and the asymptotic expansion of the solution is given by (1) with $\Theta = \pi/2$ and $\mu_l = 2(l-1)$ [1].

The subdomain $\Omega_2$ does not contain a singularity in a true sense because the derivative of the solution in the radial direction does not reach infinity when $r$ reaches 0; this can be seen by differentiating (1) with $r$ and using $\mu_k = 2(k-1)$ . However, the goal of this experiment is to analyze the behavior of the proposed method with the original computational domain decomposed into three subdomains as in Fig. 2. The parameters of this experiment that are fixed are $N = 400$ and $M = 100$, while $R_1$ and $R_2$ are varied from 0.05 to 0.5 and from 0.1 to 0.4, respectively. The above constitutes a set of experiments for a given $N_{a,\Omega_1}$ . In Fig. 10 are presented three sets of experiments for $N_{a,\Omega_1} = N_{a,\Omega_2} = 2,3,5$ with their corresponding curves clustered around the dashed curve, which in turn corresponds to the same problem setup (i.e. input parameters) but applied on the computational domain of Fig. 4 (one singularity only); it should be noted that for these set of experiments we employed constant elements instead of linear. Each cluster of curves corresponds to a different $N_{a,\Omega_1}$ and each curve from each cluster corresponds to a different $R_2$. The deviation of the solid lines from the dotted reflects the effect of the presence of $\Omega_2$ on the capacitance error of $\Gamma_6$, defined by (25). As it can be seen from Fig. 10, the effect of $\Omega_2$ for $N_{a,\Omega_1} = 2$ is minimal throughout the whole range of $R_1$ and $R_2$. For $N_{a,\Omega_1} = 3$ the effect is more prominent only at the range of $R_1$ ~0.15-0.25, while for $N_{a,\Omega_1} = 5$ the effect is prominent at the range of $R_1$ ~0.25-0.5. It is under investigation whether the discrepancies can be attributed to certain aspects of the proposed method, such as the proximity of the two singularities or the ratio of the radii of the two subdomains, when the involved aspects are quite numerous. These discrepancies may come from outside the proposed method; e.g. the accuracy of the system solver can be compromised by the condition number of the **A** matrix when $N_\alpha$ increases.

From the results of Fig. 10 we can deduce, however, that small values of $R_1$ and $R_2$ renders the proposed method less sensitive to the other input parameter ($N_\alpha$). This behavior of the proposed method is desirable in computational practice since there is no need to seek optimal values for several input parameters. Thus, the choice of small $R_1$ and $R_2$, combined with the fact that a small $R$ requires only small $N_\alpha$ (cf. Fig. 5), gives a general guideline for the preferable choice of the input parameters.

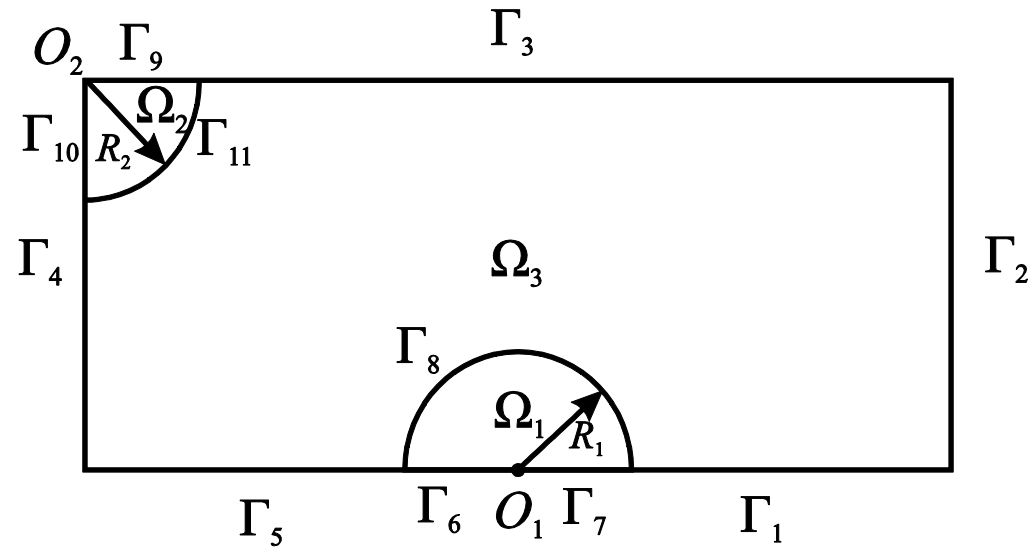

Fig. 9. Domain decomposition of the Motz problem with three subdomains.

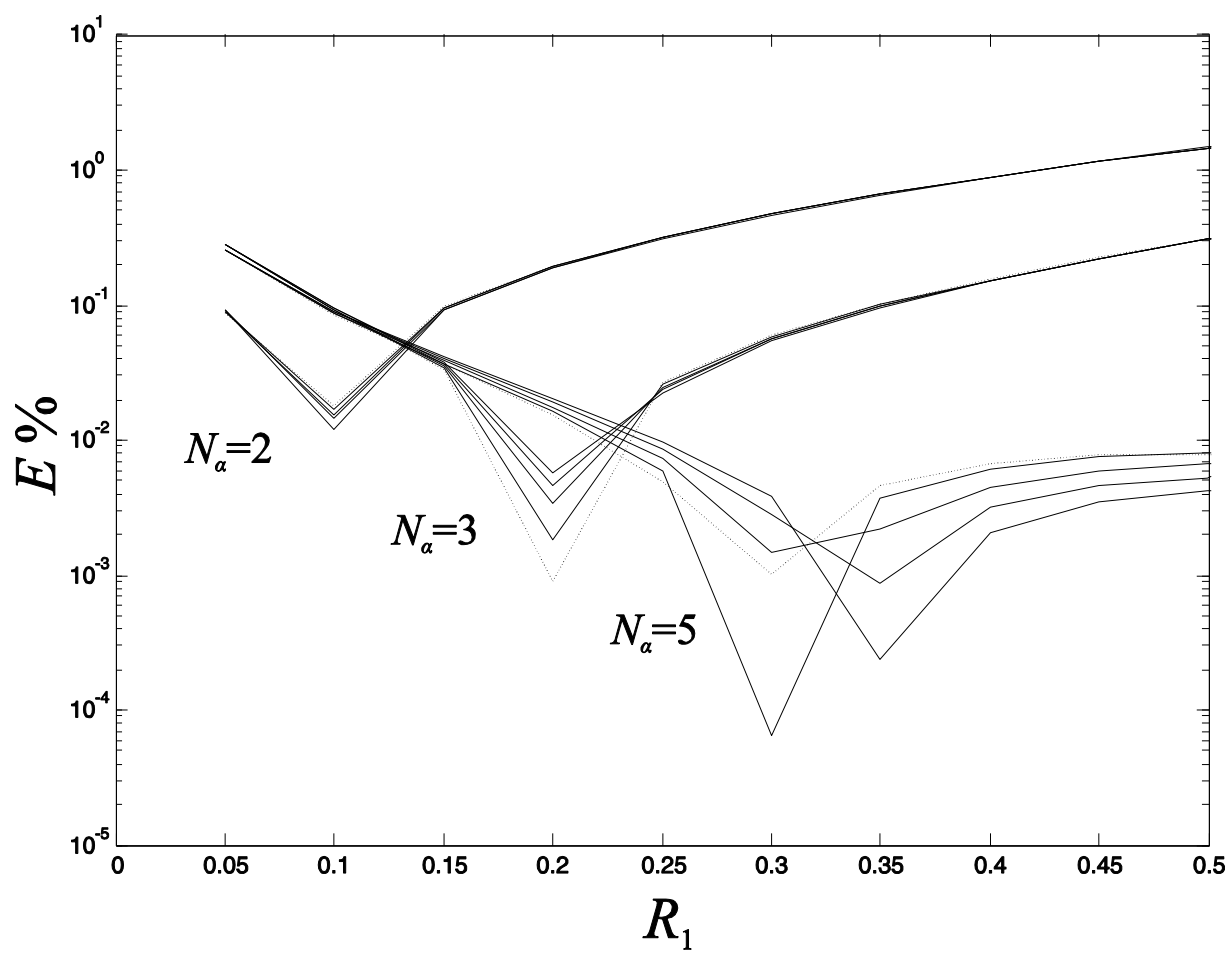

Fig. 10. Relative error, $E$, vs. $R_1$; $N$=500, $M$=100, $R_2$=0.1-0.4.

## 4 Conclusions

The hybrid BEM/SFBIM treats elliptic boundary value problems with singularities efficiently, providing adequate accuracy in the results while maintaining low computational cost in terms of domain discretization, degrees of freedom, etc. That is in contrast with the standalone BEM or GFEM, where accuracy comes at the expense of the computational efficiency. The proposed method can be seen as an augmented BEM with singular functions and it is seamlessly embedded in an existing BEM code. The method is best suited for applications that require the computation of the leading singular coefficients (generalized stress intensity factors) or the accurate computation of the derivative of the solution near singularities. The results are evaluated in terms of the capacitance of the Dirichlet boundary adjacent to the singularity and are compared with those of the BEM and GFEM.

# 5 Acknowledgements

This work has been funded by the project ΠΕΝΕΔ 2003. The project is cofinanced 80% of public expenditure through EC - European Social Fund, 20% of public expenditure through Ministry of Development - General Secretariat of Research and Technology and through private sector, under measure 8.3 of OPERATIONAL PROGRAM `COMPETITIVENESS' in the 3rd Community Support Program.